\documentclass[acmtog,nonacm]{acmart}
\setcopyright{none}
\renewcommand\footnotetextcopyrightpermission[1]{}
\usepackage[T1]{fontenc}
\usepackage[utf8]{inputenc}
\usepackage{microtype}
\usepackage[normalem]{ulem}

\usepackage{amsthm}
\usepackage{mathtools} %
\usepackage{siunitx}
\usepackage{ifthen}
\usepackage{float} %
\usepackage{hyperref}
\usepackage{cleveref}
\usepackage{csquotes}
\usepackage{calc}
\usepackage{placeins}

\usepackage{algorithmicx}
\usepackage{algpseudocode}
\usepackage{listingsutf8}

\usepackage{array}
\usepackage{booktabs}
\usepackage{multirow}

\usepackage{float}
\usepackage{wrapfig} %
\usepackage{graphicx}
\usepackage[percent]{overpic}
\usepackage{subcaption}

\usepackage{amsfonts}
\usepackage{algpseudocode,algorithm,algorithmicx}

\algrenewcommand\algorithmicrequire{\textbf{Precondition:}}
\algrenewcommand\algorithmicensure{\textbf{Postcondition:}}

\definecolor{DaviColor}{rgb}{0,0.70,0.20}
\definecolor{ZizhouColor}{rgb}{0.34,0.1,0.25}
\definecolor{DanieleColor}{rgb}{0.56,0.34,0.62}
\definecolor{DenisColor}{rgb}{0.55,0.35,0.05}

\iftrue
\newcommand{\DP}[1]{{\leavevmode\color{DanieleColor} Daniele: #1 $\qed$}}
\newcommand{\DZ}[1]{{\leavevmode\color{DenisColor} Denis: #1 $\qed$}}
\newcommand{\ZH}[1]{{\leavevmode\color{ZizhouColor} Zizhou: #1 $\qed$}}
\else
\newcommand{\DP}[1]{}
\newcommand{\DZ}[1]{}
\newcommand{\ZH}[1]{}
\fi

\newcommand{\nothing}[1]{}

\providecommand{\finalversion}{0} %
\ifthenelse{\equal{\finalversion}{1}}
{
	\renewcommand{\DP}[1]{}
	\renewcommand{\ZH}[1]{}
	\renewcommand{\DZ}[1]{}
	\renewcommand{\DT}[1]{}

}
{}

\definecolor{forestgreen}{rgb}{0.13,0.54,0.13}
\definecolor{darkblue}{rgb}{0,0,.5}
\hypersetup{
	unicode=true,
	colorlinks=true,
	citecolor=forestgreen, %
	linkcolor=darkblue, %
	urlcolor=darkblue, %
}

\renewcommand{\algorithmicrequire}{\textbf{Input:}}
\renewcommand{\algorithmicensure}{\textbf{Output:}}
\algrenewcommand{\algorithmiccomment}[1]{{\footnotesize\color{forestgreen}\(\triangleright\) #1}}

\crefname{algocf}{alg.}{algs.}
\Crefname{algocf}{Algorithm}{Algorithms}

\crefname{appsec}{Appendix}{Appendices}

\newcommand{\fprm}{q}  %
\newcommand{\prm}{\mathbf{\fprm}}  %

\newcommand{\vu}{\mathbf{u}} %
\renewcommand{\vv}{\mathbf{v}} %
\newcommand{\vx}{\mathbf{x}}  %
\newcommand{\vg}{\mathbf{g}} %

\newcommand{\per}{\theta} %

\newcommand{\vp}{\mathbf{p}}

\let\originalleft\left
\let\originalright\right
\renewcommand{\left}{\mathopen{}\mathclose\bgroup\originalleft}
\renewcommand{\right}{\aftergroup\egroup\originalright}

\renewcommand{\geq}{\geqslant}
\renewcommand{\leq}{\leqslant}

\newcommand{\N}{\mathcal{N}}

\DeclareFontFamily{U}{mathx}{\hyphenchar\font45}
\DeclareFontShape{U}{mathx}{m}{n}{<-> mathx10}{}
\DeclareSymbolFont{mathx}{U}{mathx}{m}{n}

\graphicspath{{images/}}

\begin{document}
\title{Optimized shock-protecting microstructures}

\author{Zizhou Huang}
\affiliation{%
\institution{New York University}
\country{USA}
}
\email{zizhou@nyu.edu}

\author{Daniele Panozzo}
\affiliation{%
\institution{New York University}
\country{USA}
}
\email{panozzo@nyu.edu}

\author{Denis Zorin}
\affiliation{%
\institution{New York University}
\country{USA}
}
\email{dzorin@cs.nyu.edu}

\begin{abstract}

Mechanical shock is a common occurrence in various settings, there are two different scenarios for shock protection: catastrophic protection (e.g. car collisions and falls) and routine protection (e.g. shoe soles and shock absorbers for car seats). The former protects against one-time events, the latter against periodic shocks and loads. Common shock absorbers based on plasticity and fracturing materials are suitable for the former, while our focus is on the latter, where elastic structures are useful. Improved elastic materials protecting against shock can be used in applications such as automotive suspension, furniture like sofas and mattresses, landing gear systems, etc. Materials offering optimal protection against shock have a highly non-linear elastic response: their reaction force needs to be as close as possible to constant with respect to deformation.  

In this paper, we use shape optimization and topology search to design 2D families of microstructures approximating the ideal behavior across a range of deformations, leading to superior shock protection. We present an algorithmic pipeline for the optimal design of such families combining differentiable nonlinear homogenization with self-contact and an optimization algorithm.    These advanced 2D designs can be extruded and fabricated with existing 3D printing technologies. We validate their effectiveness through experimental testing. 
\end{abstract}

\begin{teaserfigure}
	\centering
  \includegraphics[]{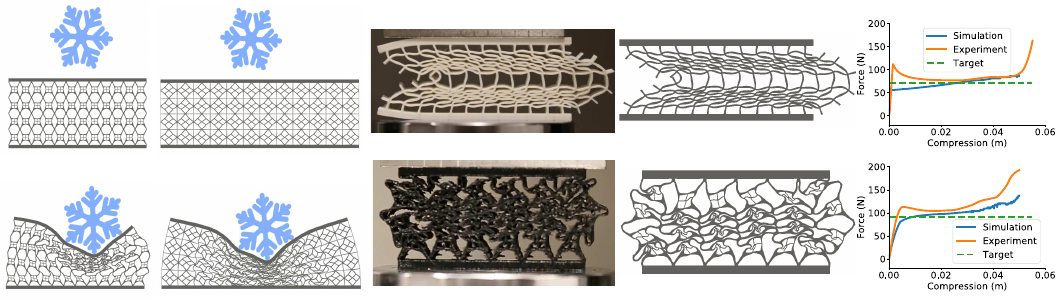}
  \caption{Our shock-protecting microstructures are designed to provide a reaction force as close as possible to constant for a wide range of displacements. We introduce a computational pipeline to design a microstructure family, and we validate its effectiveness in simulation and in physical experiments (middle). The simulated and measured response (right) is flat over a large range of compression (sample height 10cm), slowly decelerating the object and protecting it from high-impact forces.}
	\label{fig:teaser}
\end{teaserfigure}

\maketitle

\section{Introduction}
\label{sec:intro}

Mechanical shock is an abrupt and large increase in the surface force acting on an object, typically due to contact with an obstacle or another object. The need for protection from periodic mechanical shocks is common. For example, coil springs and leaf springs are extensively used in vehicle suspensions to provide a smoother ride and absorb vibrations from the road, springs can be used in robotics to act as shock absorbers, allowing robots to move more smoothly while reducing wear and tear on mechanical components, and in medical devices like prosthetics and orthotics to absorb shocks during movement. In these cases, the shock happens periodically and has a known direction.

A simple model problem, representative of most practical settings, is dropping a load with a layer of protective material on a rigid surface. The protective material layer performs two functions: first, it makes the deceleration of the object slower, reducing the maximal force acting on the load, and second, it dissipates the elastic energy to which the kinetic energy is converted.  The latter most commonly happens due to the damping/visco-elastic properties of the protective material. The former function is critical, as it eliminates the shock; the latter eliminates oscillations/prevents bounce after the initial contact. 

Minimization of maximal force/acceleration acting on the load requires materials with unusual properties: as we discuss in more detail in Section \ref{sec:problem}, the optimal behavior is for the reaction force to remain constant as the protective material deforms, which is very different from most common materials which act like springs, with the reaction force increasing with deformation. Materials with complex geometric structures, foams, or common structured materials like corrugated cardboard are commonly used as protective materials because their behavior is \emph{essentially nonlinear}\footnote{This term is commonly  used in 
the materials literature to indicate that a non-linear material model is essential to capture the material's qualitative behavior.} and closer to the ideal behavior. 

In this paper, we describe how \emph{shape optimization} for periodic microstructures consisting of 2D repeating cells can produce families of cell geometries with elastic response close to ideal over a large range of deformations, using a single base material. Our optimized structures can be fabricated using current 3D printing technologies (the 2D design is extruded to obtain a 3D shock-protecting structure), leading to shock-absorbing materials significantly closer to the perfect fixed-force deceleration. 

Our solution builds upon shape optimization algorithms for periodic metamaterials (e.g. methods producing families of cell structures spanning a particular range of effective material properties), extending them to support the distinctive features of shock-absorbing materials: 
\begin{itemize}
\item The target constitutive law is \emph{essentially nonlinear} and is not approximated well by a standard model, the homogenization must be performed in the nonlinear regime, sampling the whole stress-strain curve, rather than using a low-parametric model (e.g., captured by an elasticity tensor); 
\item \emph{Self-contact} significantly impacts the structure behavior. The deformations are large, and contact must thus be considered in the shape optimization process; 
\item Large deformations require a \emph{non-linear elasticity} model, and an accurate constitutive law for the base material must be used.
\end{itemize}

The paper shows a complete algorithmic pipeline to construct nonlinear microstructure families parametrized by the target constant stress and how to use them to realize shock-absorbing materials. 
The contributions of our paper include the following \ZH{Refactor me!}: 
\begin{itemize}
\item We formulate the equations for computing effective elastic stress-strain dependence (homogenization) of nonlinear periodic structures with cells for large displacements in the presence of contact and non-linear base material constitutive law.
\item We use a combinatorial enumeration of 2D structures to identify the best choices of structure topology for different regimes and obtain several families of cell structures with the best performance for different loads.  
\item We develop a gradient-based algorithm for shape optimization to minimize the deviation of the effective stress-strain dependence from the ideal constant-force behavior.
\item We validate the desired behavior of the resulting lattices by experimental testing of fabricated lattice samples. 
\end{itemize}

\section{Related Work}
\label{sec:related}

\paragraph{Microstructure design and optimization} There is an extensive literature on microstructure design, 
see, e.g., the survey \cite{kadic20193d} for extensive references. 

A lot of work on optimization of geometry in individual cells is based on general shape and topology optimization methods 
\cite{bendsoe1989optimal,bendsoe2003topology,allaire2002shape}. Most of these works are based on small-displacement and linear material assumptions that are fundamentally not applicable in our setting. There is an increasing number of works considering nonlinear homogenization, which we review below. 

In computer graphics, families of microstructures of various types were developed starting with \cite{schumacher2015microstructures,panetta2015elastic,Zhu:2017:TST}  with many more in studies in material science and engineering. 
We use the approach of \cite{panetta2015elastic} for our topology enumeration. Our nonlinear homogenization approach is similar to  \cite{chen2021bistable}, based on \cite{nakshatrala2013nonlinear}.

Recently, \cite{Weichen2022} used topology optimization to design microstructures to fit desired nonlinear stress-strain responses. They were also able to optimize microstructures to have a flat response. However, due to the limitation in topology optimization (Figure~\ref{fig:top-opt}) and absence of contact, they only consider moderate compressions up to 40\% and a limited range of homogenized force (10\unit{\N} to 30\unit{\N}).

\paragraph{Shock-absorbing metamaterials}
While the shock-absorbing properties of foams and structured materials were known for a long time, the desirable properties of certain types of lattices became known relatively recently. 
\cite{bunyan2019exploiting} describes the $\chi$-shaped cells,  which have flattened regions in their stress-strain curves, and which we consider as a baseline.  This type of structure was further explored in \cite{joodaky2020mechanics}. \cite{chen2020light} describes a shell-lattice metamaterial that can absorb very large energies while retaining a low density but does not attempt to optimize it. 

Early work on designing shock-absorbing structured materials \cite{kellas2010deployable} investigated designs of deployable honeycomb structures for crash energy management in light aircraft, showing these are superior to airbags. \cite{leelavanichkul2010energy} considered properties of a structure consisting of a helicoidal shell enveloping a cylinder, motivated by hydraulic shock absorbers. \cite{chen20183d} describes a new hierarchical cellular structure created by replacing cell walls in regular honeycombs with triangular lattice configurations to improve energy absorption under uniaxial compression and shape integrity at high strains.  
\cite{matlack2016composite} described elastic metastructures with wide, low-frequency band gaps while reducing global mass, with  applications in controlling structural vibrations, noise, and shock mitigation. These structures, however, are not close to the ideal shock-protecting structures we describe below. 

\cite{mueller2019energy} analyzed the energy absorption properties of various periodic metamaterials, comparing them to foam-like random structures; while random structures exhibit better uniformity of stress for varying strain, periodic lattice geometries outperform their stochastic equivalents in terms of energy absorption in some cases. We show that periodic structures can be optimized to have high stress uniformity.  More recently, \cite{acanfora2022effects} explores maximizing energy absorption in shock absorbers while minimizing thickness or mass to improve transportation safety. Their analysis is restricted to six a priori chosen structures.

\cite{GONGORA20222829} used a data-driven approach to infer the acceleration in the impact test from the stress-strain curve. The expensive transient simulation can be avoided with their approach while sacrificing some accuracy. However, they don't perform shape optimizations to find the optimal structures for impact protection.

While some works do one or two parameter sweeps to identify best-performing structures, we are not aware of any works that performed structure optimization for shock absorption systematically. 

We also briefly mention several papers that use bistable structures for shock absorption.  In this case, the transition from one stable mode to a second stable mode allows the structure to store energy and yet be reversible, assuming no plastic deformation, as pointed out in \cite{frenzel2016tailored}.  This type of structures is suitable for protection against one-time shock (e.g., a fall), but cannot protect from repeated shocks, as encountered in shipping and transportation.  
Some examples of works of this type include
\cite{izard2017optimal,ha2018design}, which describe tetra-beam-plate cells with snap-through behavior for large deflections. 
\cite{cao2021bistable} surveys a variety of bistable structures with a focus on applications to actuators, MEMS, and shock absorption. 
Most recently, \cite{jeon2022synergistic} describes a realization of a common tilted-beam bistable structure with liquid crystal elastomers (LCE), with  viscoelastic behavior improving 
energy absorption, and \cite{fancher2023dependence} proposes a biomimetic shock-absorbing mechanism inspired by the bi-stable elongation behavior of a protein.

\paragraph{Nonlinear homogenization}
Nonlinear homogenization of periodic structures for large displacements/strains is a far more complex problem than linear homogenization.  In this case, the  effective dependence between stress and strain requires multiple simulations. Even more fundamentally, for given boundary conditions for a periodic cell, the solution may be non-unique, and the material behavior may not even be fully captured by a local constitutive law.  Nevertheless, suitable approximations of effective stress-strain dependencies were obtained under certain assumptions (e.g., \cite{debotton2006neo}).
We consider a  version of the problem, with the stress-strain response for only one direction being of interest, which is considerably simpler than the general problem.  As we have mentioned above, our nonlinear homogenization approach is similar to  \cite{chen2021bistable}, based on \cite{nakshatrala2013nonlinear}, 
and used for microstructure design using topology optimization in 
\cite{wang2014design}.

We note that more general techniques for nonlinear homogenization were developed, but remain quite expensive. E.g, 
\cite{yvonnet2007reduced} and  \cite{schroder2014numerical} use reduced-order models for homogenization obtained using proper orthogonal decomposition (POD) to increase efficiency. 
These methods are further extended in \cite{fritzen2018two},  
\cite{kunc2019finite} and  \cite{kunc2020many}, with a typical approach of first constructing a reduced-order model, then sampling deformation space using this model, and finally interpolating the samples using various types of interpolation.

Several works use nonlinear, finite strain homogenization in the topology optimization context to obtain periodic metamaterials with desired properties, starting with \cite{wang2014design}, which uses numerical tests of response to a deformation, which can be considered partial homogenization, in the context of truss-based and continuum topology optimization.  A more general case of homogenization is \cite{behrou2021topology}.   \cite{wallin2020nonlinear} describes how non-linear homogenization based on  the multiscale virtual power method can be used in the context of topology optimization, with sensitivities transferred from microscale to macroscale. 
While our method is somewhat related to topology optimization methods as we use an implicit shape representation described in \cite{Panetta:2017:WCS}, unlike these techniques, it supports accurate differentiable contact.  

\cite{xue2022mapped}  uses a formulation for cellular metamaterial optimization for large deformations based on the shape map, mapping a fixed reference configuration to an optimized one.  Our method, while using an implicit shape representation, uses a similar discretization at each step to compute the shape derivatives. 

In computer graphics literature,  bistable auxetic structures are described in \cite{chen2021bistable} and used for deployable surfaces; \cite{Sperl2020} simulates yarn-level cloth effects using nonlinear homogenization; \cite{Schumacher2018} proposes a comprehensive approach to characterizing the mechanical properties of structured sheet materials with nonlinear homogenization and uses inverse design to explore structures with desired properties.

We note that a number of alternative approaches, in particular, $FE^2$ method \cite{feyel2003multilevel}  and its variations, instead of obtaining a homogenized constitutive law, use coupled two-scale simulations.  \cite{xu2022direct} merges the two simulation levels in FE$^2$ method into a single system, reducing computational complexity while maintaining accuracy.  It provides similar results to full FE meshes but with fewer degrees of freedom.

\section{Method}
\label{sec:method}

\subsection{Background and Problem formulation}
\label{sec:problem}

\begin{figure}
    \centering
    \includegraphics[width=0.4\columnwidth]{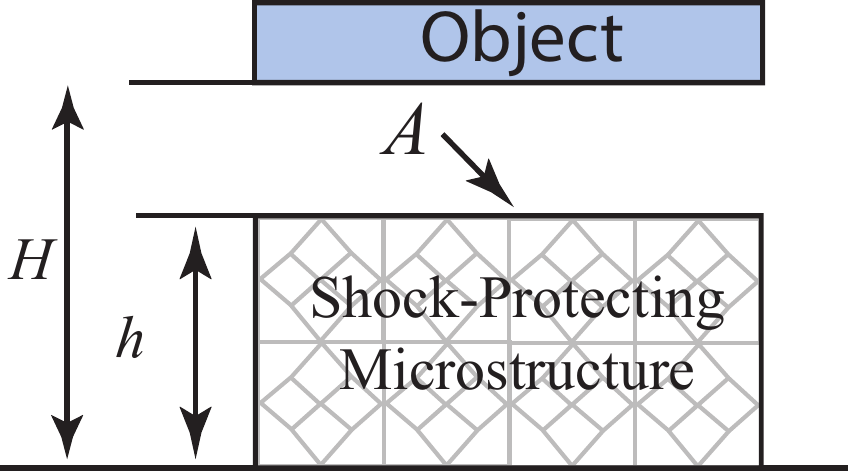}
    \caption{Model problem setup.}
    \label{fig:model-problem}
\end{figure}

We start with reviewing the problem setup (Figure~\ref{fig:model-problem}) for measuring the shock-protective properties of a material. 

A  (meta)material is typically characterized by the stress-strain curve $\sigma(\epsilon)$.  A response to a one-dimensional load is of primary interest to us, so in the model setup, we only consider one diagonal component of the stress corresponding to vertical compression, and its dependence on applied strain along the same direction i.e., we consider curve  $\sigma =\sigma(\epsilon)$, where $\sigma$ is a scalar stress, and $\epsilon$ is the scalar strain. 

\paragraph{Ideal shock-protective material.}
Suppose the kinetic energy of an object  is $m v^2/2$ right before impact, where $m$ is the object mass, $v$ is its velocity. Let  $A$ be the area of contact with the protective material.  Ignoring gravity, the force acting on the object, as the protective material is compressed to  strain $\epsilon$ is $F = A\sigma(\epsilon)$.  The assumption that the object stops for some $\epsilon < 1$, can be expressed as
\[
A h \int_0^1 \sigma(\epsilon) d\epsilon \geq mv^2/2
\]
i.e., that the work of the elastic force where $h$ is the protective layer thickness.  Here, we approximate the strain by constant over the thickness of the layer.  While in reality there may be considerable variation, this assumption is needed to obtain a problem formulation independent of the protective layer thickness/geometry. 

This leads to the following optimization problem for the "ideal"
stress-strain curve: 
\[
\min_{\sigma} \max_\epsilon \sigma(\epsilon), 
\;\mbox{subject to}\; \int_0^1 \sigma(\epsilon) d\epsilon \geq \sigma_f 
\]
where $\sigma_f = mv^2/(2Ah)$.  It is easy to see that the optimal solution is $\sigma(\epsilon)=\sigma_f$, as if  $\sigma \leq \sigma_f$ 
everywhere on $[0,1]$, the constraint can only be satisfied if the equality holds, and if $\sigma > \sigma_f$ anywhere,  this choice of $\sigma(\epsilon)$ is suboptimal, because 
the constant $\sigma$ is valid and has lower maximum. 

\paragraph{Optimization problem.}
Clearly, such a flat response is not physically possible, as close to $\epsilon=1$, the reactive forces have to increase to infinity; similarly, close to $\epsilon =0$, the reactive forces have to be close to zero.  So for any real (meta)-material, there is a ramp up part of the curve, a flat part, and a final part, corresponding to extreme compression. 

This leads to the following optimization problem: 

\emph{For a base material and a target value of stress $\sigma_f$, optimize the geometry of a unit cell so that the stress-strain curve for a metamaterial obtained by periodically repeating it is as close as possible to $\sigma(\epsilon)= \sigma_f$.} 

Solving this problem yields a one-parametric family of cell geometries $P(\sigma_f)$, parametrized by $\sigma_f$;  for each, we have an extent $\alpha(\sigma_f) < 1$ of the flat part of the stress-strain curve. 

To illustrate how such as family of materials can be used, we consider the following standard problem: given maximal allowed deceleration $G$, and expected drop height $H$, choose the optimal material in the family and required thickness.  Conservatively, assuming that all deceleration  happens at the flat part of the response curve, and approximating the strain by constant, we obtain $\sigma_f$ from the force balance $mG = A\sigma_f$; note that this does not depend on the thickness of the protective layer, and allows us to pick a material already.
For a specific material in the family corresponding  to $\sigma_f$,  we require that the work done by elastic forces on the flat part of the stress-strain curve is sufficient to absorb the kinetic energy, i.e. (approximately) 
\[
\alpha(\sigma_f) \sigma_f = \frac{mgH}{h}
\]
from which thickness $h$ can be estimated.

\subsection{Approach overview}
\label{sec:overview}

We obtain the families of protective metamaterials using a combination of combinatorial enumeration of topologies and shape optimization. 

\begin{figure}
    \centering
   \includegraphics[width=0.8\linewidth]{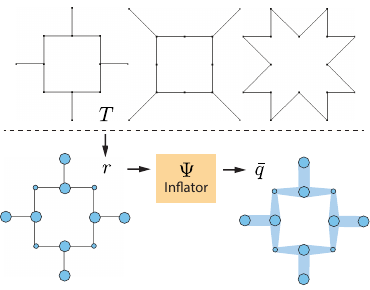}
    \caption{A cell topology $T$ is annotated with geometric parameters $r$ (a radius and a 2D position for each vertex). The inflator $\Psi$ converts the graph representation into an implicit function, which is then triangulated to obtain a mesh representation $\Bar{q}$ of the cell domain $\Omega$.}
    \label{fig:shape-param}
\end{figure}

The main components of our algorithm include: 
\begin{itemize}
    \item {\bf Topology enumeration and geometric parametrization (Section \ref{sec:inflator}).} The topology of our cells is defined by a graph within the cell, with geometric parameters are given by radii at graph nodes and blend parameters, as shown in Figure~\ref{fig:shape-param}. The outer loop of the overall algorithm enumerates different possible topologies.
    \item  {\bf Nonlinear differentiable homogenization (Section \ref{sec:hom}).}  The objective of our optimization is deviation of the stress-strain curve $\sigma(\epsilon)$ from  a constant $\sigma_f$.  To obtain 
    the effective stress $\sigma(\epsilon)$, we use periodic nonlinear homogenization with \emph{contact}, obtaining effective stresses for a set of  background deformations $\epsilon$.  Contact is of particular imporance in our setting as the material is designed specifically for very large deformations. Along with computing effective stresses, we compute their gradients with respect to shape parameters, which are essential for efficient optimization.  
    \item {\bf Objective and Optimization (Section \ref{sec:optimization}).}  For every topology and a target flat stress value $\sigma_f$, we optimize the shape to minimize the deviation of effective stress from $\sigma_f$, computed via homogenization, for strains sampled at a fixed set of strains $\epsilon_i$.  
\end{itemize}
Nonlinear differentiable homogenization is the core part of the algorithm; we discuss it in detail after briefly reviewing topology enumeratin,  starting with the forward simulation and then explaining how the derivatives can be computed. 

\subsection{Topology enumeration and geometric parametrization}
\label{sec:inflator}

We use the  graph representation  introduced in \cite{panetta2015elastic} to represent our periodic cells, and the implicit surface definition based on this graph proposed in \cite{Panetta:2017:WCS} which we briefly review here. 

\paragraph{Cell parametrization and graph inflator.}
Each microstructure is parametrized by a graph $T$, annotated with a radius and a position for every node (stacked in a single vector of parameters $r$), embedded in a rectangle of size $a \times b$; reflection symmetry is imposed on shape parameters, reducing their number.  \cite{Panetta:2017:WCS} defines an implicit surface definition that "inflates" the graph based on the radii assigned to vertices; A periodic triangular mesh of the domain $\Omega$ is obtained from the implicit function using marching squares. The map from the parameters $r$ and cell dimensions $a$, $b$ to the vector of periodic vertex positions $\bar{\prm}$  is denoted
$$\Bar{\prm} = \Psi_T(r,a,b).$$
The derivatives of the vertex positions in this mesh with respect to shape parameters (shape velocities) are computed using implicit function differentiation as explained in \cite{Panetta:2017:WCS} 
\paragraph{Design space}
The design space of the shape consists of three parts: combinatorial choice of the graph $T$, microstructure shape parameters $r$ and the size of unit cell $a$, $b$. 
Although the scale of the cell does not affect small deformation properties like the homogenized elasticity tensor, it does affect the stress-strain curve as the elasticity model is nonlinear. 

For the choice of topology, we consider 105 patterns in 2D following~\cite{Panetta2015}, generated by enumeration of patterns with bounded number over vertices in the cell and number of edges meeting at a vertex. In the optimization, we first generate the mesh in the unit square based on the graph and its parameters. Then we scale the shape by the scale parameters to get the unit cell in a rectangle.

\subsection{Nonlinear homogenization}
\label{sec:hom}

A periodic metamaterial consists of repeating cells with identical geometry $\Bar{\prm} = \Psi(r,a,b)$, 
 parametrized by the shape parameters $r$ and cell dimensions.   The effective (homogenized) properties of the material are obtained in the limit of cells repeated infinitely, and deformations are considered at a scale much larger than the cell size; in this case, we can assume that the metamaterial behaves as a homogeneous solid material, with an effective constitutive law, relating stress to strain at each point.  This stress-strain dependence can be obtained from the constitutive law of the base material and cell geometry by \emph{homogenization}.

 While for our problem the dynamic behavior of the material may be important we consider \emph{static} deformations only in our optimization, which captures most significant aspects of behavior of highly absorbing materials.  We do not include dissipation in our simulation as it does not affect the static case. We also assume negligible plasticity which is a valid assumption for materials chosen to provide protection from repeated shock.  
The deformation of the metamaterial can be decomposed into a slow-changing deformation, that can be computed from the (a priori unknown) macroscopic constitutive law, and a cell-scale fluctuation. At a level of a single cell, macroscopic stress and strain can be viewed as constant, i.e., corresponding to a linear deformation of the cell, with a periodic cell fluctuation $\Tilde{u}$ added on top.
Homogenization assumes that there is a constant macroscopic strain, equivalently, a linear deformation $Gx$ where $x$ is the spatial coordinate (test deformation), solves for the periodic fluctuation $\Tilde{u}$ and computes the resulting effective stress. 

For small displacements, a linear effective constitutive law $\sigma = C\epsilon$ can be assumed, fully determined by components of an effective elasticity tensor, which considerably simplifies the problem: the elasticity tensor can be fully inferred from a small number of test deformations. 

However, the materials we aim to construct are inherently nonlinear (Figure \ref{fig:ablation:material}).  In this case, to approximate the effective stress-strain dependence, we need to compute the effective stress resulting from a larger set of finite deformations.  The need for sampling for our problem is considerably reduced by considering only stress-strain dependence for a single direction. 

\paragraph{Notation.}   We use $x$ to denote the coordinate on the periodic cell $\Omega \subset V$, where $V$ is a rectangular tiling
with tiles of size $a \times b$. We use $u(x) = \Tilde{u}(x) + Gx$ to denote the solution of the elasticity equations with contact, 
where  $\Tilde{u}$ is the periodic fluctuation part and 
$Gx$ is the macroscopic linear deformation part, with 
$G \in\mathbb{R}^{2\times 2}$.
We restrict matrices $G$ to be symmetric to eliminate rotational components of the deformation, which do not affect elastic behavior.
The domain $\Omega$ is discretized into a periodic triangular mesh. 

The vector of coefficients of $\Tilde{u}(x)$ in a FE basis $\phi_i$ is denoted $\Tilde{\vu}$ (we use quadratic elements); this vector includes degrees of freedom only, i.e., the periodicity conditions on $\vu$ are used to exclude values on the right and upper boundaries of $V$. 

We denote the vertices of this mesh $x_i$, with the vector of vertices of size $N$ denoted $\vx$. These are determined by the shape parameters $r$ as described above. 
We use piecewise-linear basis $\xi_i$ to represent changes in the mesh as the shape parameters are varied. 
For the discrete solution $\vu$, the following equation holds: 
\begin{equation}
\vu = \Tilde{\vu} + \vx G^T
\label{eq:uc_to_u}
\end{equation}
where $\Tilde{\vu}$ is the discretization of the fluctuation. 
As explained below, the variables in the elasticity equations we solve to compute the homogenized stress-strain dependence are 
$\vu$ and components $G^{00}, G^{01}$ of the deformation matrix $G$.  
We denote the vector of all of these variables $\vv = [\Tilde{\vu}; G^{00},G^{01}]$.

\paragraph{Sampling effective stress.} As our focus is on  response to loads in a single direction, we sample a single diagonal component of the stress tensor corresponding to vertical deformation, and how it relates to the component of the macroscopic deformation $G$ in the vertical direction. 
This yields a sampled approximation of a stress-strain curve. 
In other words, the result of our homogenization procedure is a set of samples $\sigma(\epsilon_i)$ approximating the dependence 
$\sigma(\epsilon)$. 

For each sample value of $G^{11} =\epsilon$, we set up a nonlinear elasticity problem to determine corresponding stress. 
We assume that the material is free to deform in other directions; 
for this reason, we include the components $G^{00}$ and $G^{01}$ as
variables in our optimization:

\begin{equation}
 \min_{\Tilde{\vu},\ G}  W(\Tilde{\vu}, G) \quad \text{such that }  G^{11}=\epsilon 
 \label{eq:elasticity}
\end{equation}
where $W:=W_e + W_c$ is the sum of elastic (Equation \ref{eq:welastic}) and contact barrier energy (Equation 5 in \cite{Li2020IPC}). The effective stress tensor corresponding to $\epsilon$ can be computed as 
\begin{equation}\label{homogenized_stress}
\Bar{\sigma}(\epsilon) := \frac{1}{|V|}\int_\Omega \sigma(\nabla_x \Tilde{u} + G) dx,
\end{equation} 
where $|V|$ is the area of $V$ in 2D.

\paragraph{Elasticity.} The elastic energy has the form
\begin{equation}\label{eq:welastic}
W_e := \int_\Omega w_e(\nabla_x \Tilde{u} + G) dx,
\end{equation}
where $w_e: \mathbb{R}^{2\times 2} \rightarrow \mathbb{R}$ is the Neo-Hookean energy density function
$$
w_e(F - I) := \frac{\mu}{2}  (\text{Tr}[F F^T] - 2 - 2\log (\det F)) + \frac{\lambda}{2} \log^2 (\det F), 
$$
where $F$ is the deformation gradient, $\lambda$ and $\mu$ are the Lam\`e parameters. To solve \eqref{eq:elasticity} the Jacobian and Hessian of the elastic energy are needed to solve the elasticity equation. The gradient and Hessian of the elastic energy with respect to $F$ are: 
$$\sigma := \nabla w_e \in\mathbb{R}^{2\times 2}, \quad C:= \nabla^2 w_e\in\mathbb{R}^{2\times 2\times 2\times 2}.$$
Then we have the following expressions for the components of the gradient and Hessian:  
$$
\partial_{\vu_i} W_e = \int_\Omega \sigma(\nabla u):\nabla \phi_i dx\quad \partial_{\vu_i,\vu_j} W_e = \int_\Omega \nabla \phi_i : C(\nabla u) : \nabla \phi_j dx,
$$
where $i,j=1,\dots,N$.
The first and second derivatives of $W_e$ with respect to $v$ are obtained by applying the chain rule of Equation \ref{eq:uc_to_u}.

\paragraph{Periodic Contact.}
To adapt the IPC~\cite{Li2020IPC} to the periodic homogenization, we need to consider contact not only inside the periodic cell, but also between geometry from adjacent cells should be considered. 
We assume that we do not need to consider contact between geometry in cells that  are not adjacent: although theoretically this may happen, we have not observed this even for extreme deformations.
To handle contact, we use $2\times 2$ tile of the deformed periodic cell in the collision detection and barrier energy computation. 
We observe that due to periodicity, it is sufficient to consider neighbors only below and to the left of a given cell, not above and to the right.

\begin{figure}
    \centering
    \includegraphics[width=0.7\columnwidth]{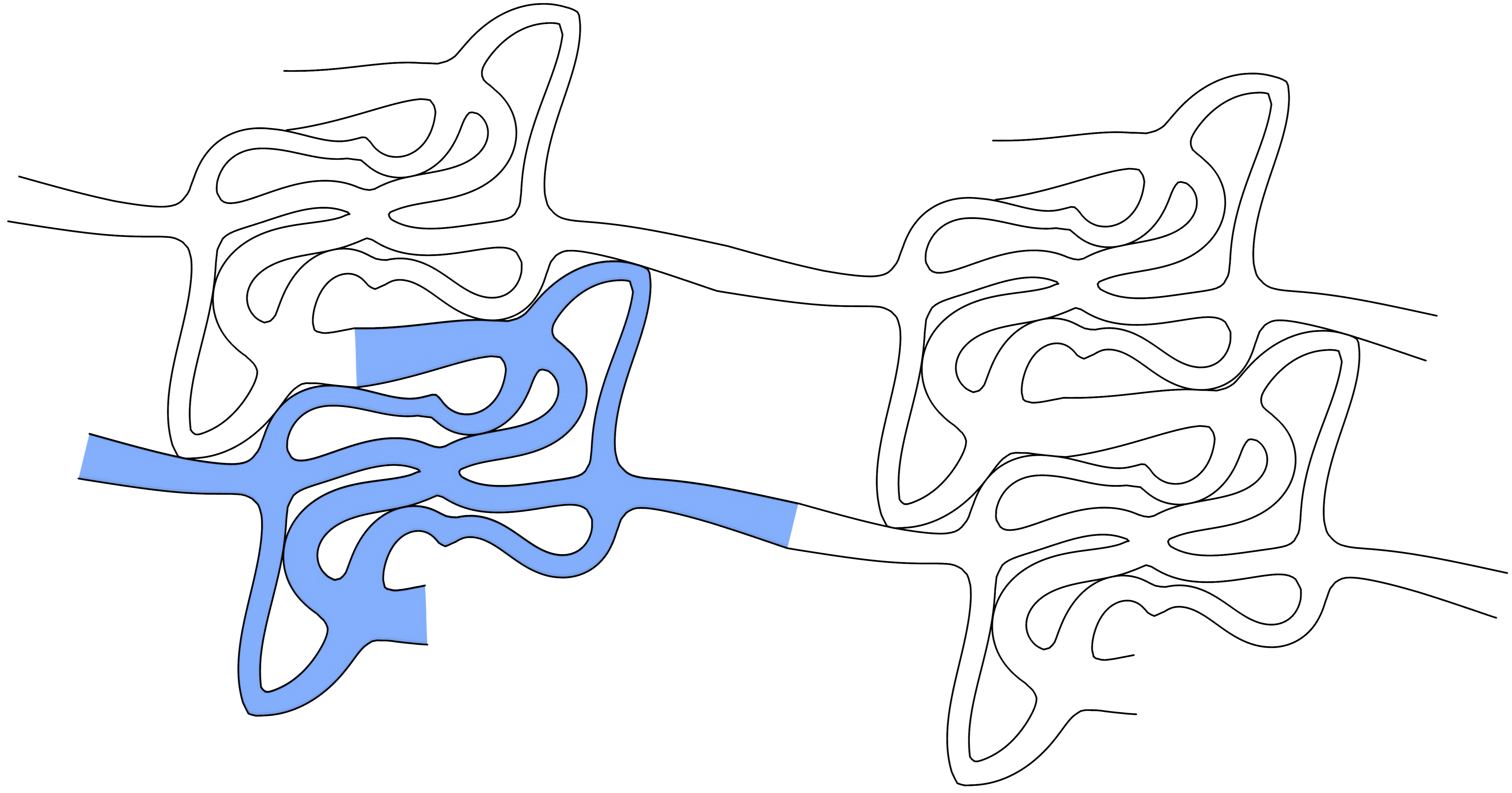}
    \caption{A deformed periodic cell collides with its tiled boundary mesh during homogenization. Accounting for collision is crucial to designing a shock-protecting microstructure family (Figure \ref{fig:ablation:contact}).}
    \label{fig:periodic-contact}
\end{figure}

Figure~\ref{fig:periodic-contact} shows how the tiled boundary mesh of a deformed periodic cell is used to detect collisions between cells.

Define $\vu^t$ to  be the vector of displacements on the $2\times 2$ tile. 
In the two-by-two tiling of copies of periodic domain $V$, 
the coordinates of vertices of three tiles are given by $\vx_i + a \mathbf{e}_1,\ \vx_i + b \mathbf{e}_2,\ \vx_i + a \mathbf{e}_1 + b \mathbf{e}_2$, where $\mathbf{e}_d$ ($d=1,2$) is the unit vector along $d$-th axis. We concatentate these along with original domain degrees of freedom into a vector $\vx^t$ of size $M$.
The index mapping $I$, maps vertex $j$ on the tiled mesh to the corresponding vertex $I(j)$ on the original mesh. The displacement on the tiled mesh $\vu^t\in \mathbb{R}^{M\times 2}$ can be represented as
\begin{equation} \label{eq:single_to_tiled}
    \vu^t_j = \Tilde{\vu}_{I(j)} + G \vx^t_j, 
\end{equation}
where $x^t_j$ is the position of vertex $j$ on the tiled mesh.

The Jacobian and Hessian of the barrier energy with respect to $u^t$ are identical to the ones used in \cite{Li2020IPC}. We apply the chain rule based on \eqref{eq:single_to_tiled} to obtain the Jacobian and Hessian with respect to $\vv$:
\begin{align*}
d_{\vv} W_c &= (d_{\vu^t} W_c) (d_{\vv} \vu^t) \\
d_{\vv}^2 W_c &= (d_{\vv} \vu^t)^T (d_{\vu^t}^2 W_c) d_{\vv} \vu^t
\end{align*}
where $d_{\vv} \vu^t$ is the gradient of the linear mapping in Equation~\ref{eq:single_to_tiled}. Entries of $d_{\vv} \vu^t$ can be computed by
\begin{align}\label{eq:single_to_tiled_grad}
    \partial_{\Tilde{\vu}_{id}} \vu^t_{jk} &= \delta_{dk} \delta_{I(j), i} \quad &i=1,\dots,N;\ j=1,\dots,M;\ d,k=1,2 \nonumber \\
    \partial_{G_{dp}} \vu^t_{jk} &= \delta_{dk} \vx_{jp}^t \quad &j=1,\dots,M;\ d,p,k=1,2
\end{align}
where $\delta_{ij}$ is the Kronecker delta.

\paragraph{Non-uniaxial Load}

\begin{figure}
    \centering
    \includegraphics[width=\linewidth]{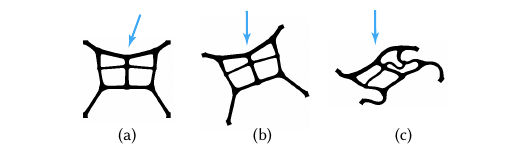}
    \caption{Instead of considering the unit cell under non-uniaxial load in (a), we rotate the rest shape and still compress along the Y direction in (b), the compressed shape is shown in (c). The load direction is shown in blue arrows.}
    \label{fig:optimize-rotated}
\end{figure}

To protect shocks in different directions, we also consider non-uniaxial loads in the nonlinear homogenization. Similar to the problem (\ref{eq:elasticity}), instead of fixing $G^{11}$, one can fix the compression strain in the load direction, which becomes a linear equality constraint. To avoid enforcing such constraint, we rotate the domain instead, so that the load is still in Y direction (Figure \ref{fig:optimize-rotated}). In this way, we don't need to change the formulation of the homogenization, only need to consider the influence of rotation in the shape derivatives. In this work, we only consider small variations around the Y direction ($\leq 20^{\circ}$) since larger angles require many more samples in strain and take a much longer time to optimize.

\subsection{Objective and optimization algorithm}
\label{sec:optimization}

\begin{figure}[h]
    \centering
    \includegraphics[width=0.7\linewidth]{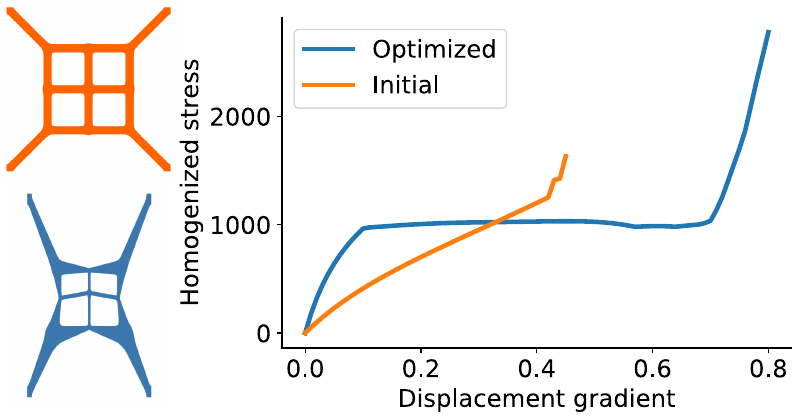}
    \caption{Each microstructure topology (orange, left) is initialized with a default set of positions and radii for each vertex. Before optimization (orange, right) the stress (\unit[per-mode = symbol]{\Pa}) - strain curve is almost linear. After optimization (blue, left) the curve is flat over a large range of deformation (blue, right).}
    \label{fig:opt}
\end{figure}

\label{sec:obj}

\paragraph{Objective}
Given a target scalar stress $\sigma_f > 0$, the goal of the optimization is to minimize the deviation of the effective stress from $\sigma_f$ on the compression strain range of $[0.1,\epsilon]$, where $\epsilon < 1$ is the max compression strain we consider. We initialize $\epsilon=0.3$ in the first optimization and gradually increase it in the following optimizations until it cannot be reached.

In every optimization, the forward simulation is solved on a series of scalar compression strains $\epsilon_i$ ($i=0,1,\dots$) uniformly sampled in the interval $[0.1,\epsilon]$, the corresponding homogenized stress and macro strain are $\Bar{\sigma}_i$ and $G_i$, the objective is 
\begin{equation}\label{eq:objective}
   J(\prm):=\sum_i (\frac{\Bar{\sigma}_i(1,1)}{\sigma^*} - 1)^2 + w (|G_i^{01}| - 0.05)_+ ^ 2, 
\end{equation}
where $\sigma^*$ is the target homogenized stress, and $\prm = (r,a,b)$ are the shape parameters and the size of the cell. The second term  penalizes shear under compression, which may lock adjacent unit cells when they shear in different directions.

Figure \ref{fig:opt} shows an example of how the stress-strain curve changes after the optimization, and how the geometry of the microstructure tile changes. 

\paragraph{Algorithm}
See Appendix \ref{app:algorithm} for the pseudocode of our complete optimization and simulation algorithms.

We use the incremental load method~\cite{OGDEN1992437} for forward simulations. However, enforcing high compression strain on an arbitrary structure may result in huge contact forces, causing convergence issues in IPC, so we optimize every structure incrementally: We first optimize the shape so that its homogenized stress reaches the target in the strain range $[10\%,30\%]$, then increase the max strain and use the previously optimized shape as the initial guess. In our experience, the simulation seldom fails in this incremental way since a structure with our target homogenized stress is unlikely to have huge contact forces inside under a slightly larger compression. 

Since we enforce reflection symmetry on the shape while the deformation is asymmetric, there exist at least two solutions (one is the reflection of the other) to the homogenization problem. This doesn't affect the convergence of our shape optimization, because the objective is the same for both solutions (homogenized stress $G^{00}$ is the same while $G^{01}$ only differs by a sign). However, the transient simulations on the microstructure tiling can be non-deterministic due to this reason, so we don't perform optimization directly on transient simulations.

\begin{figure}
    \centering
    \includegraphics[width=\linewidth]{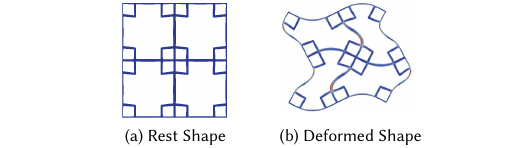}
    \caption{Homogenization on a $2\times 2$ tile of periodic unit cells. The deformed shape is not periodic in terms of every single cell.}
    \label{fig:2x2-behavior}
\end{figure}

As shown in Figure~\ref{fig:2x2-behavior}, the deformation of some microstructures in our design space cannot be fully captured by simulations on a single periodic cell, since the homogenization on a $2\times 2$ tile is not periodic in terms of every single cell. This kind of behavior also happens in other microstructures as studied in \cite{2010Bertoldi, 2022Xue}. To avoid the inconsistency caused by this behavior while keeping the runtime affordable, we perform the homogenization on the $2\times 2$ tiles after the optimization succeeds and drop shapes whose stress-strain curve is different from the curve with single cell. Then these shapes are optimized again using the more expensive $2\times 2$ homogenization. We show one example where the optimization on $2\times 2$ cells makes a difference in Section~\ref{sec:eval:ablations}.

\begin{figure}
    \centering
    \includegraphics[width=\linewidth]{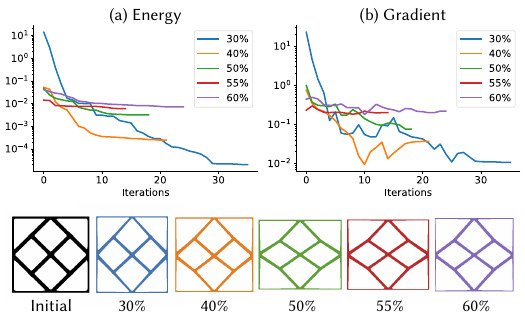}
    \caption{The objective and its gradient plots for the incremental shape optimizations (top). The initial shape and shapes after each optimization are shown (bottom).}
    \label{fig:opt-energy}
\end{figure}

Figure~\ref{fig:opt-energy} shows an example of how the objective and its gradient reduce in the shape optimization. Since our goal is to reduce the objective until the homogenized stress is close enough to the target, for better efficiency, we stop the optimization when the point-wise error is small even if the gradient is not small enough.

\subsection{Shape derivatives}
\label{sec:shape-deriv}

We optimize the objective \eqref{eq:objective} of Section~\ref{sec:obj} (in the pseudocode $\Call{Objective}{\prm, \Tilde{\vu}, G}$) with respect to the microstructure cell parameters $\prm = (r,a,b)$. In the following, we will overload the notation (whenever the derivation is generic) and use $\prm$ to refer to either the refers to the independent vertex positions $\Bar{\prm} = \Psi(r)$ defining domain $\Omega$ (we eliminate a subset of boundary vertices due to periodicity) or to the cell dimensions $\Tilde{\prm} = (a,b)$.
As the cell scale is typically determined by fabrication constraints we can use box constraints to keep one of the scale dimensions close to the desired size. 
Shape parameters (graph vertex positions and radii) determine vertex positions, as described in Section~\ref{sec:inflator}.

The derivatives of the objective $J$ with respect to the shape parameters are computed using the adjoint method, with which the shape derivatives can be obtained by solving a single additional linear equation, and then evaluating an expression depending on this unknown. See Appendix~\ref{app:shape-derivatives} for the complete derivation of shape derivatives.

\section{Evaluation}
\label{sec:results}

We implemented our algorithm in C++ and used Eigen~\cite{eigenweb} for the linear algebra routines, a modified version of PolyFEM~\cite{polyfem} for finite element simulation,  triangle \cite{shewchuk2005triangle} for meshing, and Pardiso~\cite{pardiso-7.2a,pardiso-7.2b,pardiso-7.2c} for solving linear systems. All our experiments are run on a cluster node with an Intel Cascade Lake Platinum 8268 processor limited to use 16 threads.

We first show the coverage of our microstructure family in the space of (strain, stress) pairs: a point $(\sigma_f, \epsilon)$ is considered covered if for strain $\epsilon$ the actual response of the microstructure in the family corresponding to $\sigma_f$ does not deviate from $\sigma_f$ by more than 10\%. We show representative examples of shock-protecting lattices, which are fabricated using a Prusa Mk3s printer in TPU (sample size 10cm tall, 2.6cm thick) and physically tested under compression using a INSTRON 5966 Mechanical universal testing machine (Section \ref{sec:eval:family}). We limit the physical validation to a subset of our microstructure topology due to the time required for each test (\textasciitilde24 hours printing time per sample). %
We then provide a comparison against the closest microstructure found in previous work (Section \ref{sec:eval:baseline}), and conclude the evaluation with ablations (Section \ref{sec:eval:ablations}) for the use of a non-linear material model, of a contact model, and for restricting the homogenization to a single axis.

\subsection{Microstructure Family}
\label{sec:eval:family}

\begin{figure}
    \centering
    \includegraphics[width=0.95\linewidth]{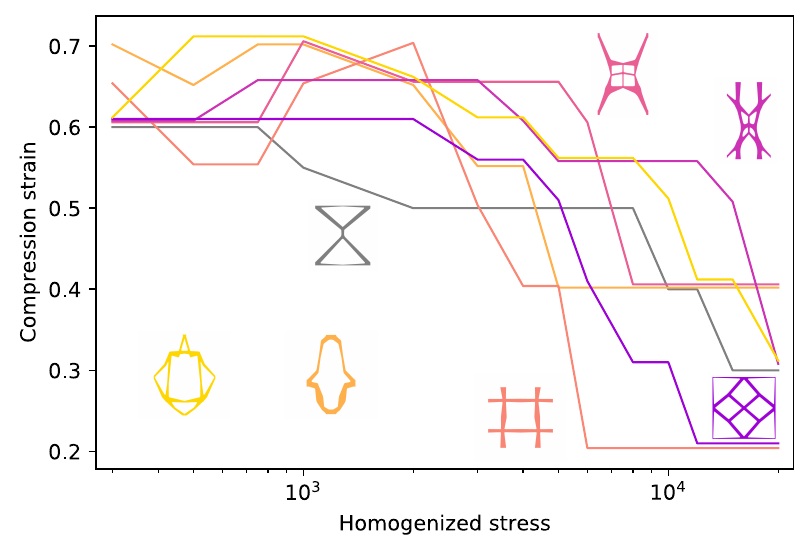}
    \caption{For each target stress value, we plot the maximal possible strain for which the structures' homogenized stress equals to the target. The 6 shown topologies (and their coverage) are selected from a total of 105 topologies to maximize their coverage, while providing a compact representation. We include the baseline structure from \cite{joodaky2020mechanics} (grey) to show that our family has much wider coverage than the state-of-the-art.}
    \label{fig:family}
\end{figure}

\begin{figure}
    \centering
    \includegraphics[width=0.9\linewidth]{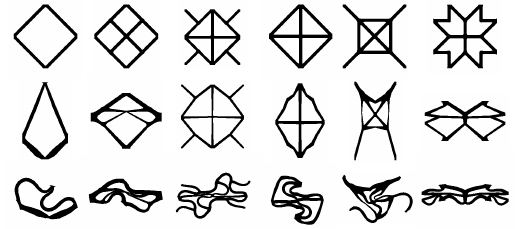}
    \caption{Examples of rest initial shapes (top), rest optimized shapes (middle), and compression of optimized shapes (bottom).}
    \label{fig:more-opt-patterns}
\end{figure}

\begin{figure}
    \centering
    \includegraphics[width=\linewidth]{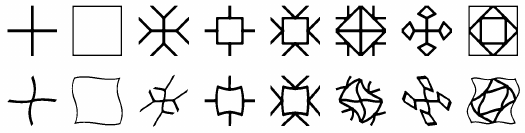}
    \caption{Topologies failed to fit a constant stress-strain curve up to 30\% strain. Initial shapes (top) and compression of initial shapes (bottom) are shown.}
    \label{fig:failed-patterns}
\end{figure}

To find the material coverage of a microstructure topology, we select 14 homogenized stress targets (from 300 to 20000) and run our incremental optimization to find parameters for a flat response curve for 7 different  compressive deformations (from 20\% to 70\%). We ran this procedure for 105 topologies, which took approximately 145k CPU hours of computation. Then we filtered the curve to find a subset of 6 providing a good coverage (Figure \ref{fig:family}). More examples of optimized shapes are shown in Figure~\ref{fig:more-opt-patterns}. Some of the failed topologies are shown in Figure~\ref{fig:failed-patterns}. %

\paragraph{Compression Tests.}

We validate our microstructure family by creating a rectangular object tiled with a small number of cells and performing a compression test. We perform the compression test virtually (creating a single triangular mesh of the entire object, and simulating it using PolyFEM, using a Neo-Hookean material model, backward Euler time integration, and with contact) and physically, using a universal testing machine.

\begin{figure*}
    \centering
    \includegraphics[width=0.96\linewidth]{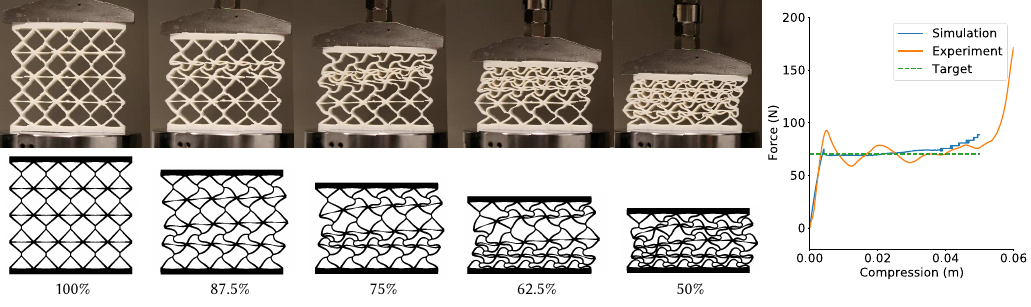}
    \caption{The microstructure topology proposed in \cite{joodaky2020mechanics} is fabricated after our optimization. Due to complex buckling, the simulation (bottom left) and experiment (top left) don't match exactly, but both the deformed shapes and the force-displacement plots (right) are similar. The thickness of the compressed shape is marked as a percentage with respect to the initial thickness.}
    \label{fig:compression}
\end{figure*}

Figure \ref{fig:teaser} shows a representative example of three families in our coverage: the virtual compression tests show a good agreement with our homogenized target, and the physical experiments confirm that our physical models are correctly modeling the real-world deformation of an isotropic material (we used thermoplastic polyurethane for this experiment). In Figure \ref{fig:compression}, we show the same experiment for the extended $\chi$-shaped structure proposed by \cite{joodaky2020mechanics}, optimized with our algorithm. 

We note that the flat region of force-displacement curves in the compression tests is not as wide as in the periodic homogenization, because in most examples the top and bottom half rows of the microstructure tiling are not periodic, so not able to deform as the periodic cells. The flat region can be widened by stacking more rows, but our fabrication is limited by the size of our 3D printer.

Although the simulation results are similar to the experiment results, they don't match exactly for the following reasons: Since the unit cell is symmetric while its deformation is asymmetric, the solution is not unique (one solution is the reflection of the other), so both experiment and simulation results can not be uniquely determined by the boundary conditions and may be influenced by small perturbations, e.g. anisotropy caused by 3D printing, floating point errors and non-determinism caused by parallelism in simulations.

\begin{figure}
    \centering
    \includegraphics[width=\linewidth]{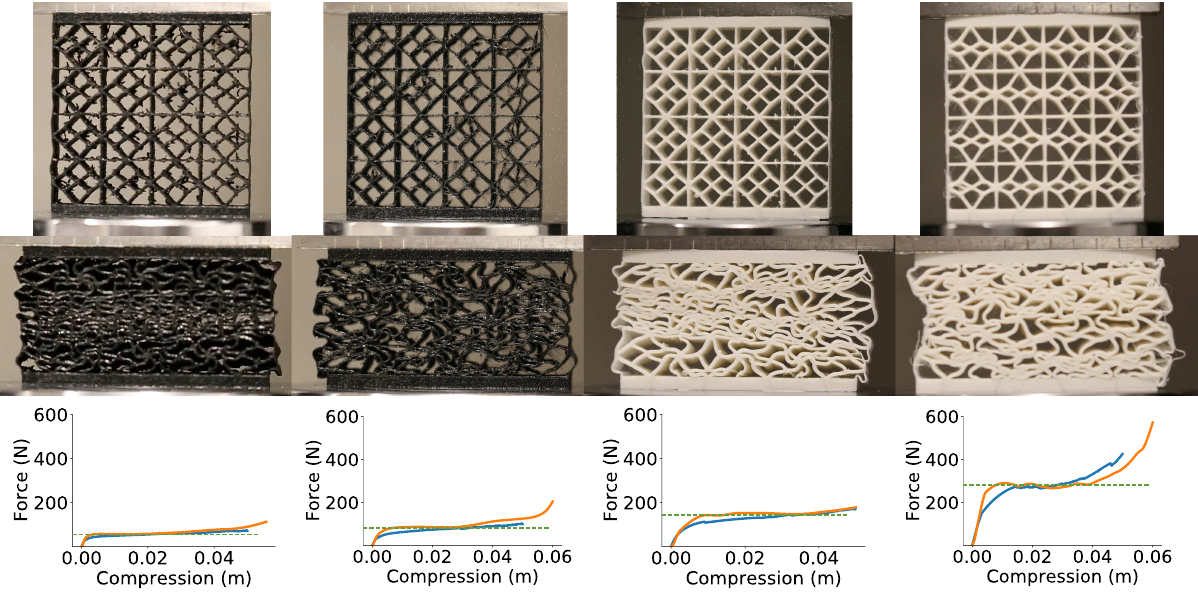}
    \caption{The Rhomboids unit cell optimized for different effective stresses. We show the rest shape (top row), the geometry after applying a compression of 50 mm (middle), and the force - compression distance curves of experimental data (orange), simulation (blue), and target force (green).}
    \label{fig:compression-single}
\end{figure}

The geometric variations within a single family are subtle (Figure \ref{fig:compression-single}), but lead to very different response curves. Physical validation results are in line with our computational predictions, with a close match on the response curve, despite pushing the resolution of our 3D printer to the limit (many features in our printed sample use only one or two lines of plastic due to resolution limitations of FDM printing).

\begin{figure}
    \centering
    \includegraphics[width=0.96\linewidth]{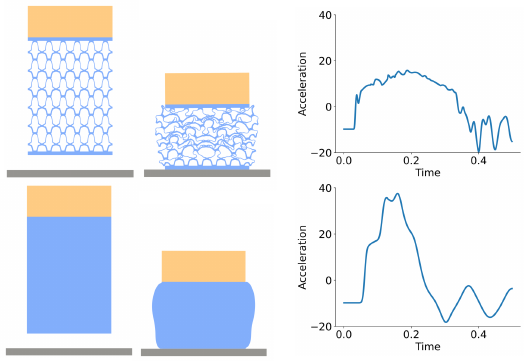}
    \caption{Drop test with the microstructure stitched to the bottom of the object, simulated on one of our optimized structures (top) and a solid material of the same dimension with smaller Young's modulus (bottom) so that they compress by the same distance at highest compression. We show the plots of time (\unit{\second}) and acceleration (\unit[per-mode = symbol]{\meter\per\square\second}) on the right, where we can observe that our maximal acceleration is around half of a solid block of a soft material.}
    \label{fig:drop}
\end{figure}

\subsection{Drop Tests.}
\paragraph{Baseline Comparison}
To evaluate the use of our approach to design protective gears or packaging protection, we run simulated drop test experiments, where the microstructure is attached to the falling object (Figure \ref{fig:drop}). We observe that the flat response of our structure leads to slower deceleration of the load. In contrast to shock protectors relying on plastic material~\cite{acanfora2022effects}, our structure returns to its rest state after impact, making it reusable.

\begin{figure}
    \centering
    \includegraphics[width=\linewidth]{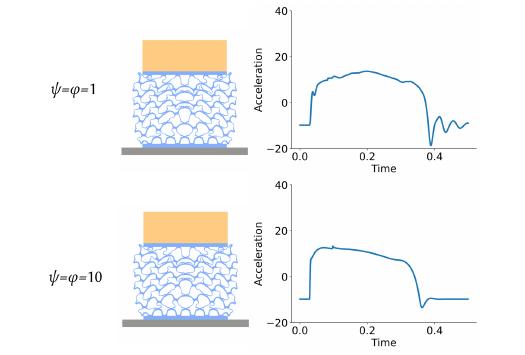}
    \caption{Drop test (see Figure~\ref{fig:drop}) with damping.}
    \label{fig:drop-damp}
\end{figure}

\paragraph{Visco-elasticity}
To analyze the influence due to visco-elasticity, we add the strain-rate proportional damping from ~\cite{damping2018} into our transient simulations. The maximum acceleration with damping is very similar to the value without damping (Figure~\ref{fig:drop-damp}), but the oscillations of the curves are reduced. However, optimizing the transient simulations with visco-elasticity, though feasible, is much more unstable and expensive; for this reason, following previous works, we do not include the dynamics effects, and as a consequence, visco-elasticity in the optimization. 

\begin{figure}
    \centering
    \includegraphics[width=\linewidth]{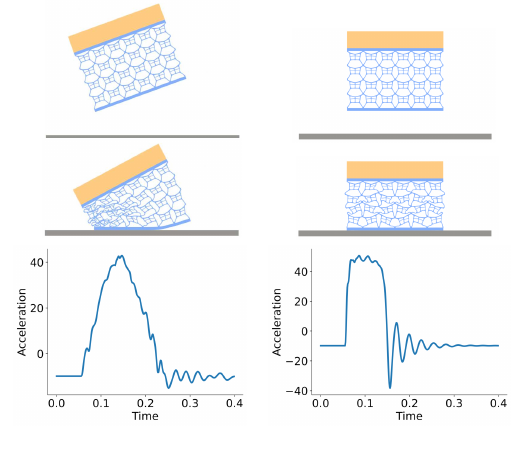}
    \caption{Drop tests with a rotated load (left) and a horizontal load (right) with the same density. We show the plots of time (\unit{\second}) and acceleration (\unit[per-mode = symbol]{\meter\per\square\second}) at the bottom.}
    \label{fig:drop-rotated}
\end{figure}

\paragraph{Rotated Load}
To evaluate the effectiveness of the optimized structures under shocks in perturbed directions, we pick one topology, perturb the load direction by $20^{\circ}$, and optimize the stress-strain curves to be constant in both the vertical direction and the perturbed direction. We then simulate the drop test (Figure \ref{fig:drop-rotated}) with a vertical load and a load rotated by $20^{\circ}$, with the same density (Figure \ref{fig:drop}). The maximal acceleration is similar in two scenes.

\subsection{Package Protection}
\label{sec:eval:package}

\begin{figure}
    \centering
    \includegraphics[width=\linewidth]{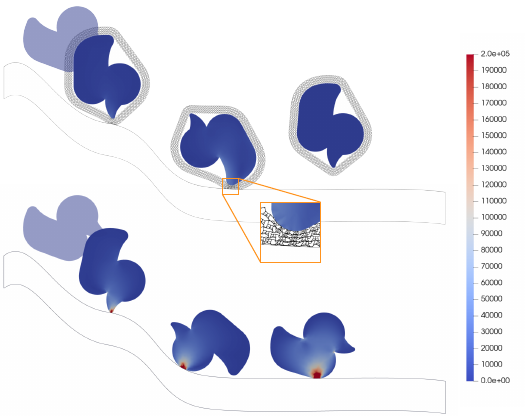}
    \caption{A duck falls down the slope with and without our protecting microstructures. The Von Mises stress distribution is shown.}
    \label{fig:duck}
\end{figure}

\begin{figure}
    \centering
    \includegraphics[width=\linewidth]{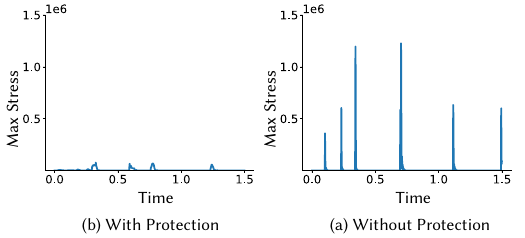}
    \caption{Max Von Mises stress plots of Figure~\ref{fig:duck}. The maximum stress is reduced by 94\% by adding our protecting microstructures.}
    \label{fig:duck-plot}
\end{figure}

This example shows how our optimized microstructures can be applied in package protection. We pick the optimized shape in Figure~\ref{fig:drop-rotated}, which can protect from shocks in perturbed directions. We then generate a quadrilateral mesh covering the object and map the unit cell to every quadrilateral to form a protective shell (Figure~\ref{fig:duck}). As shown in Figure~\ref{fig:duck-plot}, under the protection of microstructures, the maximum stress on the object is reduced from \qty{1.2e6}{\pascal} to \qty{7.6e4}{\pascal}.

\subsection{Baseline Comparisons}
\label{sec:eval:baseline}
We compare with the state-of-the-art structure proposed in \cite{joodaky2020mechanics}. This extended $\chi$-shaped structure has been discovered by manual design of its topology, and by grid searching over a low parametric parametrization of its shape.

\begin{figure}
    \centering
    \includegraphics[width=0.9\linewidth]{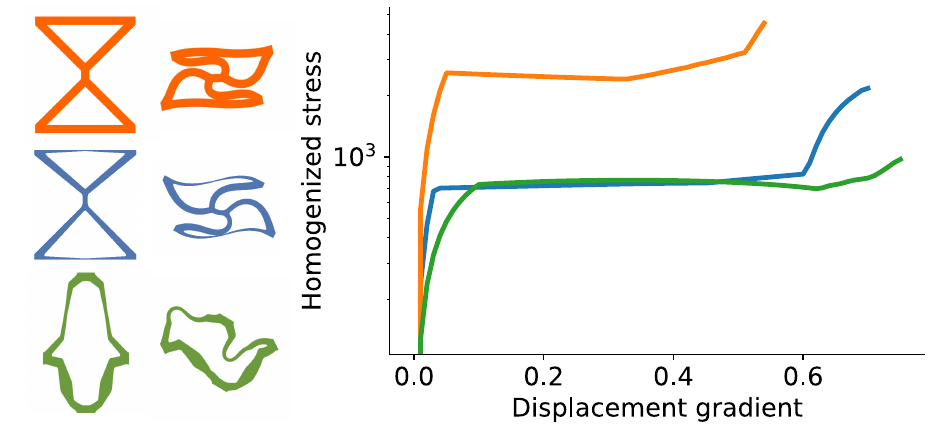}
    \caption{The flat response of the topology and geometry proposed in \cite{joodaky2020mechanics} (orange), can be considerably extended by optimizing its geometric parameters using our approach (blue). Our microstructure family has a different connectivity (green), which provides an even wider flat response. We show on the left the geometry of the corresponding cells in rest pose and compressed at 50\%.}
    \label{fig:baseline}
\end{figure}

In Figure \ref{fig:baseline}, we show that our optimization approach can modify its effective stress value by optimizing geometric parameters. Additionally, our extensive search of 105 topologies led to connectivities that can achieve higher compression with a flat response, extending the 55\% of the baseline up to 70\% (Figure \ref{fig:baseline} and \ref{fig:family}).

\subsection{Drop Experiment}
\label{sec:eval:drop-exp}

\begin{figure}
    \centering
    \includegraphics[width=\linewidth]{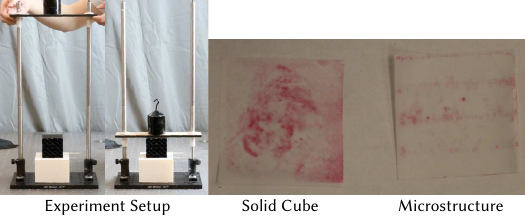}
    \caption{Drop experiment setup (left) and pressure sensing films after the drop (right). Dark color indicates high stress.}
    \label{fig:drop-exp}
\end{figure}

To further validate the shock protecting effect of our microstructure family, we perform a drop experiment on one of our optimized microstructures (Figure~\ref{fig:drop-exp}). It's hard to accurately measure the maximum force during the drop since the collision normally happens in very short time. As a result, we use the Fujifilm pressure indicating film to measure the maximum stress distribution on the weight instead. We fabricate the microstructure in Figure~\ref{fig:compression-single} of size \numproduct{6 x 6 x 6}\ \unit{\centi\meter} with black TPU and attach a piece of the pressure indicating film (with measure range 2 to 6 \unit{\kilo\gram\per\square\centi\meter}) to its bottom. After the drop, we qualitatively indicate the maximum stress based on the darkness of the color, and compare the result with the baseline, which replaces the microstructure with a 100\% filled cube made of the same material. The result shows that our microstructure significantly reduces the maximum stress during the drop.

\subsection{Ablations}
\label{sec:eval:ablations}

We provide a series of ablation experiments to motivate our choice of shape representation, including contact forces in homogenization, using a non-linear material model, and using $2\times 2$ tile in homogenizations.

\paragraph{Shape Representation}

\begin{figure}
    \centering
    \includegraphics[width=\linewidth]{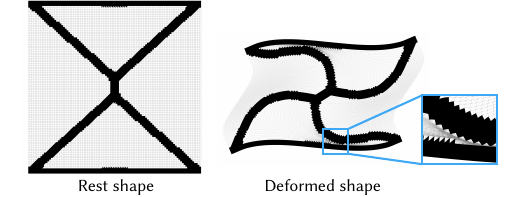}
    \caption{Homogenization simulated with topology optimization representation.}
    \label{fig:top-opt}
\end{figure}

Topology optimization is widely used in shape design for microstructures \cite{Weichen2022, WALLIN2020103324, ZHANG2019490}. To justify our choice of shape representation (Figure~\ref{fig:shape-param}) in the optimization, we perform the nonlinear homogenization on the optimized extended $\chi$-shaped structure using topology optimization. We rasterize the shape and assign Young's modulus $E=10^6$\unit{\Pa} to solid cells and $E=10^{-2}$\unit{\Pa} to void cells. In Figure~\ref{fig:top-opt}, the shape of void cells is close to singular when solid cells approach contact, resulting in convergence issues and poor accuracy in homogenization.

Although there are recent works \cite{contact-top-opt} able to resolve contact in topology optimization in some cases, robust and accurate handling of contact with topology optimization is largely an open problem as contact behavior can be significantly altered e.g. by artifacts of surface extraction.

\begin{figure}
    \centering
    \includegraphics[width=0.96\linewidth]{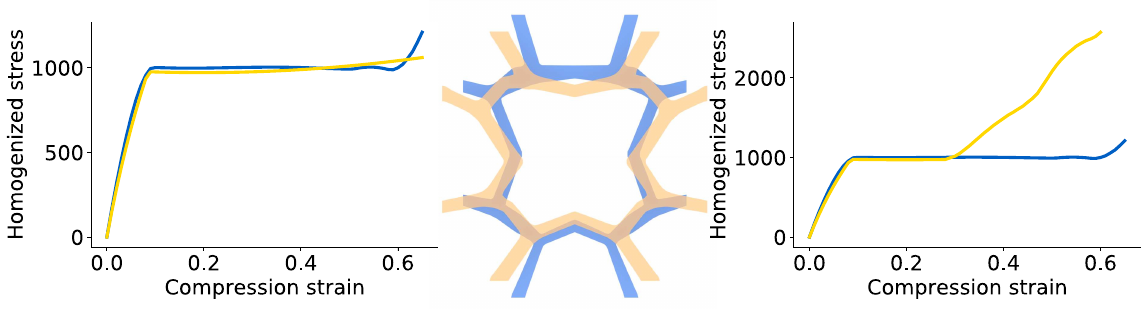}
    \caption{The same microstructure (middle) is optimized to have a flat response (plot on the left is the stress-strain plot after optimization) with (blue) and without contact (yellow). On the right, we show the stress-strain plot of the optimized structures with contact enabled. Ignoring contact during optimization leads to considerably worse performance, as taking contact into account is essential for correctly handling large deformations.}
    \label{fig:ablation:contact}
\end{figure}

\paragraph{Contact}
We compared the stress-strain plot of a pattern optimized with and without contact forces (Figure \ref{fig:ablation:contact}). While the optimization succeeds in both cases (and even makes more progress without contact forces), the contact-aware homogenized stress of the shape optimized without contact is much higher than expected. The contact forces unavoidably introduced in the test lead to a non-flat response for the material optimized without contact, while they are flat for the specimen optimized with contact.

\begin{figure}
    \centering
    \includegraphics[width=\linewidth]{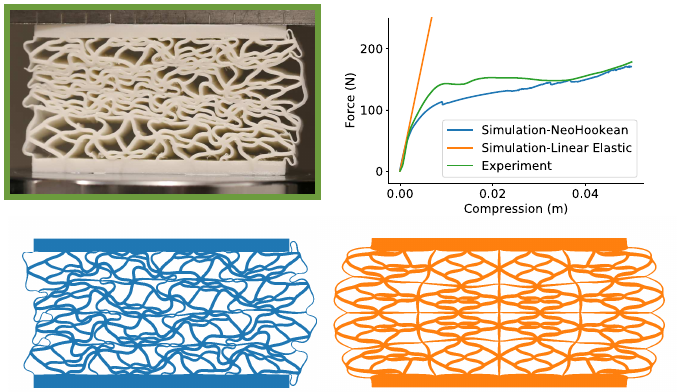}
    
    \caption{Physical experiment on a $4\times 4$ tile of the Rhomboids unit cell. Top right: force - compression distance curves of experiment data (green), simulation with NeoHookean material (blue), and simulation with Linear Elasticity (orange). Note that the simulation with Neo-Hookean material (bottom left) closely matches the compression experiment (top left), while the simulation with linear elasticity (bottom right) fails to capture the deformation.}
    \label{fig:ablation:material}
\end{figure}

\paragraph{Material Model}
We advocate for using a non-linear constitutive law for the base material to more accurately capture large deformations. In Figure \ref{fig:ablation:material} we reproduce the physical testing of one of our samples using PolyFEM configured with the linear elastic model and the non-linear Neo-Hookean model, using the same material parameters. For low compression rates, both models are accurate. After ~5\% compression the linear model diverges from the experimental data, while the non-linear model closely matches the measured force-compression curve, confirming that a non-linear material model is essential for shock-protecting materials.

\section{Conclusions}
\label{sec:concl}

We presented a method for shape optimization of homogenized microstructure materials accounting for large deformation and contact. We demonstrated its effectiveness in designing a family of shock-absorbing microstructures and validated it with simulated and physical experiments.

Our work opens the door to optimizing metamaterial families with non-linear materials and contact forces, and there are many exciting directions for future works: (1) extension to shapes with complex geometrical boundaries, for example using rhombic cells \cite{tozoni2021optimizing}, (2) support multi-axial loads, (3) account for dynamic inertial effects and/or friction in the homogenization, (4) add a plasticity model, and (5) extend the construction to 3D microstructures.

There are several limitations to our method. First, since the homogenization problem is nonconvex and highly nonlinear, multiple saddle points and local minima exist, which sometimes cause line search failure in shape optimization. Second, to obtain such a constant stress-strain response, the range of the homogenized force is limited to two orders of magnitude (Figure ~\ref{fig:family}). Third, some compressed microstructure is not completely periodic (Figure~\ref{fig:drop}), resulting in a loss of efficiency in protection.

To foster reproducibility and adoption of this technique, we will release both our microstructure family and a reference implementation of our optimization pipeline.



\appendix
\section{Shape Derivatives}
\label{app:shape-derivatives}

We will first summarize the adjoint method applied to differentiating $J$, then the computation of required partial derivatives.

\paragraph{Adjoint Method.}  

At the minima of  $W$, in Equation \ref{eq:elasticity}, we have $\partial_{\vv} W = 0$. Differentiating both sides of the equation with respect to the shape parameters $\prm$ we obtain
$$
d_{\prm} \partial_{\vv} W + \partial^2_{\vv} W d_{\prm} \vv = 0,
$$
i.e.
$$
d_{\prm}\vv = -(\partial^2_{\vv} W)^{-1}\partial_{\prm} \partial_{\vv} W,
$$
where $\partial^2_{\vv} W$ is the hessian of the total energy in the forward simulation.

The gradient of $J$ with respect to $\prm$ then can be written as
\begin{equation}
\begin{split}
    d_{\prm} J &= \partial_{\prm} J + \partial_{\vv} J d_\prm \vv \\
    &= \partial_{\prm} J -\partial_{\vv} J\ (\partial^2_{\vv} W)^{-1}\partial_{\prm} \partial_{\vv} W.
\end{split}
    \label{eq:Jderiv}
\end{equation}

Suppose $p$ is the solution of the following linear equation
$$
p^T \partial^2_{\vv} W = -\partial_{\vv} J,
$$
then it directly follows from \eqref{eq:Jderiv}, $d_\prm J$ can be simplified to
$$
d_{\prm} J = \partial_{\prm} J + p^T \partial_{\prm} \partial_{\vv} W.
$$
Thus, to compute the shape derivative we need to compute $\partial_{\vv} J$, $\partial_{\prm} J$, $\partial_{\prm} \partial_{\vv} W$ and solve a linear equation with the same coefficient matrix as the linear system  in the forward newton solve.

\paragraph{Shape derivative of elastic force $\partial_{\prm} \partial_{\vv} W_e$.}
Since the elastic energy is in the form of an integral over the domain:
$$
W_e = \int_\Omega w_e(\nabla \Tilde{u} + G) dx,
$$
we can first compute the shape derivatives with respect to all vertices in $\vx$, then apply the chain rule of $\prm\rightarrow \vx$, 
\begin{align*}
\partial_{\vx} \partial_{\vv} W_e &= \partial_{\vx} (\partial_{\vu} W_e \partial_{\vv} \vu) \\
&= (\partial_{\vx} \partial_{\vu} W_e) \partial_{\vv} \vu + \partial_{\vu} W_e (\partial_{\vx} \partial_{\vv} \vu)
\end{align*}
where $\partial_{\vx} \partial_{\vv} \vu$ can be obtained by differentiating Equation \ref{eq:uc_to_u}, and for $\partial_{\vx} \partial_{\vu} W_e$ we follow the derivation in ~\cite{huang2022differentiable}, here we only write down the final formula. Recall that $\xi_i$ is the linear  basis used to represent th change of the rest mesh, then
\begin{align*}
\partial_{\vx_i} \partial_{u_j} W_e &= \partial_{\vx_i} \int_\Omega \sigma : \nabla \phi_j dx \\
&=\int_\Omega -\sigma \nabla \xi_i^T : \nabla \phi_j - \nabla \phi_j : C : \nabla u \nabla \xi_i + (\sigma : \nabla \phi_j )  \nabla \cdot \xi_i
\end{align*}

\paragraph{Shape derivative of contact force $\partial_{\prm} \partial_{\vv} W_c$.}
The periodic contact force can be considered as a function of vertices of the tiled mesh, we follow ~\cite{huang2022differentiable} to compute the derivatives of the contact force with respect to vertex positions on the tiled mesh, then apply the chain rule to the map $\prm\rightarrow \vx^t$, which we discuss below.

For every vertex $j$ on the tiled mesh, its position can be written as
\begin{equation}\label{eq:periodic_to_tiled}
\vx^t_j = \Bar{\prm}_{I(j)} + \Tilde{\prm}^T \begin{pmatrix} \alpha_j & 0 \\ 0 & \beta_j \end{pmatrix}\quad \text{for some $\alpha_j, \beta_j\in \{0, 1\}$}.
\end{equation}
where $(\alpha_j,\beta_j)$ are indices of the tile to which the vertex belongs to, i.e., shifts by $a$ or $b$ the components of the scale part of $\prm$. To recall, $I(j)$ maps the index of vertices on the tiled mesh back to the index of vertices on the single cell mesh.
Then 
\begin{align*}
\partial_{\Bar{\prm}_i} \partial_{\vv} W_c &= \partial_{\Bar{\prm}_i} (\partial_{\vu^t} W_c\ \partial_{\vv} \vu^t) \\
&= \partial_{\vu^t} W_c\ \partial_{\Bar{\prm}_i} \partial_{\vv} \vu^t + \partial_{\Bar{\prm}_i} \partial_{\vu^t} W_c\ \partial_{\vv} \vu^t \\
&= \partial_{\vu^t} W_c\ \partial_{\Bar{\prm}_i} \partial_{\vv} \vu^t +  (\sum_{I(j) = i} \partial_{x^t_j} \partial_{\vu^t} W_c)\partial_{\vv} \vu^t\\
\partial_{\Tilde{\prm}} \partial_{\vv} W_c &= \partial_{\Tilde{\prm}} (\partial_{\vu^t} W_c\ \partial_{\vv} \vu^t) \\
&= \partial_{\vu^t} W_c\ \partial_{\Tilde{\prm}} \partial_{\vv} \vu^t + \partial_{\Tilde{\prm}} \partial_{\vu^t} W_c\ \partial_{\vv} \vu^t \\
&= \partial_{\vu^t} W_c\ \partial_{\Tilde{\prm}} \partial_{\vv} \vu^t + (\sum_j \begin{pmatrix} \alpha_j & 0 \\ 0 & \beta_j \end{pmatrix} \partial_{\vx^t_j} \partial_{\vv} W_c) \partial_{\vv} \vu^t
\end{align*}
The gradient with respect to $\prm$ can then be obtained by applying the chain rule to the map $\prm \rightarrow [\Bar{\prm}, \Tilde{\prm}]$. In the above equations, the terms we did not  dicuss yet are 
$\partial_{\Bar{\prm}_i} \partial_{\vv} \vu^t$ and $\partial_{\Tilde{\prm}} \partial_{\vv} \vu^t$, which can be computed by combining Equations \eqref{eq:single_to_tiled_grad} and \eqref{eq:periodic_to_tiled}. 

\paragraph{Derivatives of effective stress $\partial_{\prm} J,\partial_{\vv} J$.} 
To compute derivatives of Equation \ref{eq:objective} with respect to $\prm$ and $\vv$, the only difficulty is in 
$$
\Bar{\sigma}^{11}=\frac{1}{|V|}\int_\Omega \sigma^{11}(\nabla \Tilde{u} + G) dx.
$$
Similar to the elastic force, it is also in the form of an integral of the stress tensor over the unit cell domain, except that the gradient of basis function is replaced by identity, so one can derive derivatives following the derivation for elastic forces.
\section{Algorithm}
\label{app:algorithm}

The main function is {\sc Optimization}$(\prm_0,\epsilon)$, where $\prm_0$ are the initial shape parameters, and $\epsilon$ is the list of values of vertical strain for which we evaluate the stress-strain curve. After every optimization finishes, we plot the stress-strain curve with dense samples (every 1\%) to verify the optimized result: if the homogenized stress at every sample point in the range is within $10\%$ of the target stress, we accept the optimized result; otherwise, we reject it and stop optimizing for larger strain range.
 
 Function $\sc{Solve}$ solves a sequence of problems for increasing deformations $G^{11}$, using the previous result as initialization, and calling {\sc IncrementalSolve}, which imposes a constraint on $G^{11}$ as a penalty with increasing weight, as this leads to a more reliable optimization behavior. {\sc NewtonSolve} is a standard Newton method, with line search ensuring that there are no self-interesections or element inversion \cite{Li2020IPC}. It uses the {\sc ForceSPD} function to ensure that the Hessian approximation used in the solve for the descent direction is always positive-definite. \textsc{ConstrainedNewtonSolve} is similar to \textsc{NewtonSolve}, but with $G^{11}$ fixed to the input scalar strain $\epsilon$.

The algorithm uses a few auxiliary functions for which we do not provide explicit pseudocode as they are either standard or described in other papers:
\begin{itemize}
     \item \textsc{Inflate} is the mapping from shape parameters to the discretized domain (Section 6 of \cite{Panetta:2017:WCS}); 
     \item \textsc{LBFGSB} returns the descent direction using the L-BFGS solver \cite{wieschollek2016cppoptimizationlibrary} with box constraints;
     \item \textsc{LineSearch} is the standard back-tracking line search.
\end{itemize}

 \begin{algorithm}
    \begin{algorithmic}
        \Function{Optimization}{$\prm_0,\ \epsilon $} 
            \State $\prm \gets \prm_0$ 
            \State $I \gets 0$ \Comment{Number of iterations}
            \State $\Omega \gets \Call{Inflate}{\prm}$ \Comment{Section \ref{sec:inflator}}
            \State $\Tilde{\vu}, G \gets \Call{Solve}{\Omega, \epsilon}$ \Comment{Section \ref{sec:hom}}
            \State $J \gets \Call{Objective}{\prm, \Tilde{\vu}, G}$ \Comment{Section \ref{sec:optimization}}
            \State $\vg \gets \nabla \Call{Objective}{\prm, \Tilde{\vu}, G}$ \Comment{Section \ref{sec:shape-deriv}}
            \Repeat
                \State $\vp \gets \Call{LBFGSB}{J, \vg}$ \Comment{Descent Direction}
                \State $\alpha \gets \Call{LineSearch}{\prm, \vp}$
                \State $\prm \gets \prm + \alpha \vp$
                \State $\Omega \gets \Call{Inflate}{\prm}$ \Comment{Section \ref{sec:inflator}}
                \State $\Tilde{\vu}, G \gets \Call{Solve}{\Omega, \epsilon}$ \Comment{Section \ref{sec:hom}}
                \State $J \gets \Call{Objective}{\prm, \Tilde{\vu}, G}$ \Comment{Section \ref{sec:optimization}}
                \State $\vg \gets \nabla \Call{Objective}{\prm, \Tilde{\vu}, G}$ \Comment{Section \ref{sec:shape-deriv}}
                \State $I \gets I + 1$
            \Until{$\|\vg\| < \varepsilon_0$ or $\|J\| < \varepsilon_1$ or $I > IterMax$}
            \State \textbf{return} $q$
        \EndFunction
    \end{algorithmic}
\end{algorithm}

\begin{algorithm}
    \begin{algorithmic}
        \Function{Solve}{$\Omega, \epsilon$}
            \State $\Tilde{\vu}[0] \gets 0$
            \State $G[0] \gets 0$
            \For{$k \gets 1$ to $length(\epsilon)$}
                \State $\Tilde{\vu}[k], G[k] = \Call{IncrementalSolve}{\Tilde{\vu}[k-1], G[k-1], \epsilon[k]}$
            \EndFor
            \State \textbf{return} $\Tilde{\vu}, G$ \Comment{List of solutions}
        \EndFunction
    \end{algorithmic}
\end{algorithm}

\begin{algorithm}
    \begin{algorithmic}
        \Function{IncrementalSolve}{$\Tilde{\vu}_0,\ G_0,\ \epsilon_\text{Target}$}
            \State $\Tilde{\vu} \gets \Tilde{\vu}_0$
            \State $G \gets G_0$
            \State $w \gets w_0$ \Comment{Initial weight, initially $10^4$}
            \State $e_0 \gets |G^{11} - \epsilon_\text{Target}|$
            \State $e \gets e_0$
            \Repeat 
                
                \State $\Tilde{\vu}, G \gets \Call{NewtonSolve}{\epsilon_\text{Target}, w, \Tilde{\vu}, G}$
                \State $e \gets |G^{11} - \epsilon_\text{Target}|$
                \State $w \gets 2w$

                \If{$e > e_0$} \Comment{Worse than initial solution}
                    \State $w_0 \gets w$
                    \State $[\Tilde{\vu}, G] \gets [\Tilde{\vu}_0, G_0]$
                \EndIf
            \Until{$e < \epsilon_2$}

            \State $G^{11} \gets \epsilon_\text{Target}$
            \State $\Tilde{\vu}, G \gets \Call{ConstrainedNewtonSolve}{\epsilon_\text{Target}, \Tilde{\vu}, G}$ \Comment{Fix $G^{11}$ in the solve}
            \State \textbf{return} $\Tilde{\vu}, G$
        \EndFunction
    \end{algorithmic}
\end{algorithm}

\begin{algorithm}
    \begin{algorithmic}
        \Function{NewtonSolve}{$\epsilon,\ w,\ \Tilde{\vu},\ G$}
            \State $\Tilde{W} \gets W(\Tilde{\vu}, G) + w |G^{11} - \epsilon|^2$ \Comment{Energy to be minimized}
            \Repeat
                \State $H \gets \Call{ForceSPD}{\nabla^2 \Tilde{W}(\Tilde{\vu}, G)}$
                \State $\vp \gets -H^{-1} \nabla \Tilde{W}(\Tilde{\vu}, G)$ \Comment{Descent Direction}
                \State $\alpha \gets \Call{ConstrainedLineSearch}{\Tilde{\vu}, G, \vp}$
                \State $[\Tilde{\vu}, G] \gets [\Tilde{\vu}, G] + \alpha \vp$
            \Until{$\|\nabla W(\Tilde{\vu}, G)\| < \varepsilon_1$}
            \State \textbf{return} $\Tilde{\vu}, G$
        \EndFunction
    \end{algorithmic}
\end{algorithm}

\begin{algorithm}
    \begin{algorithmic}
        \Function{ForceSPD}{$H$}
            \State $\beta \gets 10^{-8}$
            \State $\Tilde{H} \gets H$
            \While{$\Tilde{H}$ is not SPD}
                \State $\Tilde{H} \gets H + \beta I$
                \State $\beta = 2\beta$
            \EndWhile
            \State \textbf{return} $\Tilde{H}$
        \EndFunction
    \end{algorithmic}
\end{algorithm}

\end{document}